# Typical periods or typical time steps? A multi-model analysis to determine the optimal temporal aggregation for energy system models

Maximilian Hoffmann[a,b,1], Jan Priesmann[c], Lars Nolting[c], Aaron Praktiknjo[c], Leander Kotzur[a], Detlef Stolten[a,b]

[a] Institute of Energy and Climate Research, Techno-economic Systems Analysis (IEK-3), Forschungszentrum Jülich, 52428 Jülich, Germany

[b] Chair for Fuel Cells, RWTH Aachen University, c/o Institute of Electrochemical Process Engineering (IEK-3), Forschungszentrum Jülich GmbH, Wilhelm-Johnen-Str., 52428 Jülich, Germany

[c] RWTH Aachen University, Institute for Future Energy Consumer Needs and Behavior (FCN), Mathieustr. 10, 52074 Aachen, Germany

## Abstract

Energy system models are challenged by the need for high temporal and spatial resolutions in order to appropriately depict the increasing share of intermittent renewable energy sources, storage technologies, and the growing interconnectivity across energy sectors.

This study evaluates methods for maintaining the computational viability of these models by analyzing different temporal aggregation techniques that reduce the number of time steps in their input time series. Two commonly-employed approaches are the representation of time series by a subset of single (*typical*) time steps, or by groups of consecutive time steps (*typical periods*). We test these techniques for two different energy system models that are implemented using the Framework for Integrated Energy System Assessment (FINE) by benchmarking the optimization results based on aggregation to those of the fully resolved models and investigating whether the optimal aggregation method can, *a priori*, be determined based on the clustering indicators.

The results reveal that *typical time steps* consistently outperform *typical days* with respect to clustering indicators, but do not lead to more accurate optimization results when applied to a model that takes numerous storage technologies into account. Although both aggregation techniques are capable of coupling the aggregated time steps, *typical days* offer more options to depict storage operations, whereas *typical time steps* are more effective for models that neglect time-linking constraints.

In summary, the adequate choice of aggregation methods strongly depends on the mathematical structure of the considered energy system optimization model, and *a priori* decisions of a sufficient temporal aggregation are only possible with good knowledge of this mathematical structure.

[1] Corresponding author: Email address max.hoffmann@fz-juelich.de

# 1. Introduction

The implementation of low-carbon technologies in the energy sector is an important aspect of the global endeavor to decelerate climate change. This includes the capacity expansion of intermittent renewable energy sources, the integration of storage technologies, as well as the realization of synergies across sectors known as 'sector coupling' [1, 2]. Consequently, the complexity of energy systems has been steadily growing, which also requires more advanced energy system models. On the other hand, the growth rate of transistor density, known as Moore's law [3, 4], has declined in recent years [5]. Additionally, the growing energy demand of high-performance computing counteracts its benefits for achieving a carbon-neutral and energy-efficient future.

Therefore, a multitude of techniques to reduce the computational complexity of energy system models have emerged [6, 7]. Temporal aggregation takes a prominent role among these, which strive to reduce the amount of input data from time series without significantly affecting the solutions of energy system models. With respect to temporal aggregation techniques found in the literature, two frequently utilized approaches can be differentiated: The first is the representation of the entire time series by single *typical time steps*, which is often referred to as *system states* or *snapshots*, and the other one is the representation of time series via groups of consecutive time steps, so-called *typical periods*. In cases where *typical periods* comprise a duration of 24 h, they are also referred to as *typical days* or *representative days*. Both approaches, *typical time steps* and *typical periods,* reduce the number of time steps that must be analyzed by energy system models.

Temporal resolution is especially important for bottom-up energy system optimization models (ESOMs), which optimize the operation of single components in different locations and their interaction with each other, and are often capable of considering the potential capacity expansion of certain technologies [8, 9]. ESOMs, which exclusively focus on the system operation of predefined technologies, are defined as dispatch models, whereas models that also optimize the design of energy systems are called capacity or generation expansion planning models [10]. Furthermore, ESOMs differ in terms of their spatial coverage, sectoral coverage, and representation of storage technologies.

The multitude of different ESOMs on the one hand and different methods available to perform temporal aggregation on the other raises the following research questions:

- Can an optimal choice of temporal aggregation technique be made *a priori* based on its capacity to represent the original time series?
- Are certain aggregation techniques predetermined for certain types of ESOMs and underlying research questions?

In an attempt to answer these, we construct herein a comparative framework that enables us to fix all influences apart from the type of ESOM and the temporal aggregation method. Under these *ceteris paribus* conditions, we hold the modeling framework, aggregation algorithms, and computing resources for performing the model runs constant. We analyze two fundamentally different models, namely a single-node model for expansion planning and a multi-node dispatch model, and made use of the computing resources provided by the *Jülich Supercomputing Centre* to perform a large-scale analysis to investigate the interdependencies between model settings and optimal aggregation methods.

The remeinder of this paper is structured as follows: In Section 2, the stated aggregation techniques of interest are reviewed with respect to their applications in the literature. In Section 3, the models are introduced. In Section 4, we summarize and analyze the results of applying both aggregation techniques to both models. Finally, our conclusions are presented in Section 5.

# 2. Review of temporal aggregation approaches for energy system optimization models

This section provides an overview of some of the temporal aggregation techniques found in the literature with respect to the models for which they are used. Here, two basic ideas must be distinguished: The representation of the original time series by a subset of single time steps, called *system states* or snapshots (hereinafter referred to as *typical time steps*); and the representation by periods of interconnected time steps (hereinafter referred to as *typical periods*), which have strengths and weaknesses. Although the use of *typical time steps* enables the diversity of time steps to be represented in a comparably small number of states, *typical periods* require more time steps, as they consist of multiple consecutive time teps but therefore enable: (1) to depict chronologies within the periods; and (2) to more accurately represent short-term dynamics. Table 1 lists some of these models and their areas of implementation for temporal aggregation.

***Table 1. Literature overview of models using temporal aggregation for energy system optimization.***

| Reference | Reference | Year | TA method | TA specification | Model type | Spatial resolution | Sectoral focus |
|---|---|---|---|---|---|---|---|
| **Monts et al.** | **[11]** | 1991 | Sys. St. | 3 Sys. St. per season | No model (just for generating marginal costs) | No model (just for generating marginal costs) | No model (just for generating marginal costs) |
| **van der Weijde et al.** | **[12]** | 2012 | Sys. St. | 500 Sys. St. | Expansion | Multi-node | Single-commodity |
| **De Sisternes Jimenez and Webster** | **[13]** | 2013 | Typ. Per. | 4 Typ. Weeks | Dispatch | Single-node | Single-commodity |
| **Wogrin et al.** | **[14]** | 2014 | Sys. St. | 6 Sys. St. | Dispatch | Single-node | Single-commodity |
| **Agapoff et al.** | **[15]** | 2015 | Sys. St. | 5 to 100 Sys. St. | Expansion | Multi-node | Single-commodity |
| **Fitiwi et al.** | **[16]** | 2015 | Sys. St. | 15 to 310 Sys. St. | Dispatch | Multi-node | Single-commodity |
| **Munoz et al.** | **[17]** | 2015 | Typ. Per. | 1 to 300 Typ. Days | Dispatch | Single-node | Single-commodity |
| **Lythcke-Jørgensen et al.** | **[18]** | 2016 | Sys. St. | 49 Sys. St. | Dispatch | Single-node | Multi-commodity |
| **Nahmmacher et al.** | **[19]** | 2016 | Typ. Per. | 1 – 100 Typ. Days | Expansion | Multi-node | Single-commodity |
| **Ploussard et al.** | **[20]** | 2016 | Sys. St. | 10 to 400 Sys. St. | Dispatch | Multi-node | Single-commodity |
| **Wogrin et al.** | **[21]** | 2016 | Sys. St. | 98 Sys. St. | Dispatch | Single-node | Single-commodity |
| **Bahl et al.** | **[22]** | 2017 | Sys. St. | 1 to 12 Sys. St. | Expansion | Single-node | Multi-commodity |
| **Härtel et al.** | **[23]** | 2017 | Sys. St. | 68, 137, 274, 548, 1095, 2190, 4380 Sys. St. | Expansion | Multi-node | Single-commodity |
| **Bahl et al.** | **[24]** | 2018 | Typ. Per. | Iterative increase of the number of typical days and time steps per day | Expansion | Single-node | Multi-commodity |
| **Tejada-Arango et al.** | **[25]** | 2018 | Sys. St. | 100 Sys. St. | Dispatch | Multi-node | Single-commodity |
| **Hilbers et al.** | **[26]** | 2019 | Sys. St. | 480, 1920, 8760 Sys. St., 60 most expensive subsamples | Expansion | Single-node | Single-commodity |
| **Priesman et al.** | **[6]** | 2019 | Res. red. | 1 hr, 2hrs, 4hrs | Dispatch | Multi-node | Single-commodity |
| **Savvidis et al.** | **[27]** | 2019 | Res. red. | 15 min, 1 hr | Dispatch | Single-node | Single-commodity |
| **Tejada-Arango et al.** | **[28]** | 2019 | Typ. Per. | 1, 2, 4 Typ. Days per month<br>1 48-h period per month<br>1 96-h period per month | Dispatch | Single-node | Multi-commodity (water and electricity) |

## 2.1. Using typical time steps

*Typical time steps* constitute a subset of time steps within a time series. They are either randomly chosen from the original time series or selected using clustering techniques. Pöstges et al. [29] demonstrated that, for very simple economic dispatch models based on peak load pricing, the input time series can be reduced to a subset of time steps without changing the built capacities and total costs of the optimal energy system.

Dispatch models, which also incorporate technologically highly-resolved unit commitment (UC) models, focus on the operation of a given system with fixed capacities. They represent technical restrictions with a high level of detail in order to be able to, e.g., accurately simulate the dispatch behavior of generation units or resulting market prices. Some dispatch models refrain from aggregating time series at all. In such cases, the marginal utility of additional accuracy always exceeds the marginal utility of additional complexity. As many dispatch models analyze narrow time horizons and incorporate temporally volatile data, many studies have been conducted on the influence of the temporal resolution within these narrow time horizons, such as by O'Dwyer et al. [30], Pandžić et al. [31], Deane et al. [32], Savvidis et al. [27], and Priesmann et al. [6]. However, *typical time steps* are also used in dispatch models that address the full time horizon. Lythcke-Jørgensen et al. [18] proposed a method for identifying representative time slices independent of their time of occurrence. Bahl et al. [22] aggregate the input time series into time slices using k-means clustering.

Apart from dispatch models that require storage behavior to be appropriately depicted, several authors have employed *typical time steps* for testing predefined systems under heterogeneous system conditions, e.g., Fitiwi et al. [16], Munoz et al. [17], Ploussard et al. [20], and Tejada-Arango et al. [25] applied the system state approach to different IEEE reliability systems [33, 34]. In the latter case, the sequence of *typical time steps* was also considered with respect to model storage units as well.

Other publications applied this reduction technique to models of the British electricity system [12, 26], electricity models for the generation expansion [15] of transmission expansion [23], and other single-node electricity models [14, 21].

In contrast, only a small number of publications have selected this approach for multi-commodity CHP systems [18, 35].

## 2.2. Using typical periods

A second way to reduce the number of time steps is the representation of the time series by a set of periods consisting of interconnected time steps. One of our proceeding studies [10] showed that the most common period lengths are days, as many input time series feature a daily pattern, but other period lengths, such as weeks [13, 36-42], are also utilized.

Temporal aggregation into *typical periods* is less frequently applied in dispatch contexts than in expansion models. This is partly due to the relatively short time horizons that most dispatch models cover (e.g., days instead of years in the case of expansion models). If long time horizons are being analyzed, modelers tend to draw on *typical periods* to reduce the complexity of the model instead of using *typical time steps*. In an early attempt, Monts et al. [11] analyzed the impact of different representative periods on the operation and resulting prices of a generation portfolio in Austin, Texas. *Typical days* were used as *typical periods* for the ELMOD model [43] and by Nahmmacher et al. [19]. De Sisternes Jimenez and Webster [13] clustered the time series into typical weeks to reduce the computational complexity of the unit commitment problem. Tejada-Arango et al. [28] combined *typical time steps* with *typical periods* (days and weeks) in a dispatch model, considering a water basin with three different reservoirs.

Temporal aggregation based on *typical periods* is often used for optimization models of residential or district energy systems or CHP systems containing multiple commodities and conversion units. The publications to be found in the literature that apply this approach are numerous and a detailed overview of them is provided by Hoffmann et al. [10]. To highlight a few examples, Yokoyama et al. [44, 45] selected a subset of *typical days* from three seasons (winter, summer, and midseason) for building models. Likewise, Wakui et al. [46-48] chose typical season days to build models with thermal and battery storage units. In turn, Mavrotas et al. [49] chose either one average day for each month or season days with a lower temporal resolution for the modeling of a hospital in Greece. Lozano et al. used monthly average days for another hospital model [50] or district model [51], in the former case with a further distinction made between weekdays and weekend days. Similarly, Harb et al. [52] selected monthly average days for a building and district model and took thermal storage units into account. The differentiation between weekdays and weekend days was also incorporated by Samsatli et al. [53, 54] and is often used when the input time series also has a strong weekly pattern, i.e., in the case of residential electricity demand. Casisi et al. [55], Weber et al. [56], Mehleri et al. [57, 58], and Wouters et al. [59] utilized three typical season days for a district model with multiple commodities, in some cases with daily storage units [56, 58, 59]. Finally, Fazlollahi et al. [60, 61] applied k-means clustering to determine *typical days* for district heating models.

More recently, clustering techniques to determine *typical periods* have gained popularity, as they decrease the deviation from the original time series compared to heuristic approaches. Marquant et al. [62, 63], Schiefelbein et al. [64], Bahl et al. [24], and Baumgärtner et al. [65] chose k-medoids clustering for district heating [62-64] and CHP [24] models. Lin et al. [66], Gabrielli et al. [67, 68], and Baumgärtner et al. [69] employed k-means clustering in cogeneration models for buildings [66] and single-node district CHP systems with different daily and seasonal storage units [67-69]. Schütz et al. [70], Kotzur et al. [36], and Zatti et al. [71] compared different clustering techniques for a residential energy system [70], three different hypothetical energy systems [36], and a district system for a university campus [71], respectively. Scott et al. [72] used a PAM algorithm to determine typical days for a national capacity expansion model. Nahmmacher et al. [19], Welder et al. [73], and Kannengießer et al. [74] used hierarchical clustering for multi-nodal energy systems containing intraday storage units and, in the latter two cases, seasonal storage technologies as well.

It is worth noting that the temporal resolution within *typical periods* can also be decreased to further reduce the total number of time steps [10]. This was performed, for instance, by Mavrotas et al. [49], based on the gradients between time steps, clustering techniques as demonstrated by Fazlollahi et al. [60], or modified PAM algorithms [24, 65, 69]. Other approaches were taken by Weber et al. [56] and Mehleri et al. [57] to reduce intradaily temporal resolution, whereas Raventos et al. [75] compared the impact of clustering typical days to clustering adjacent time steps.

Apart from that, several national multi-node electricity models have been developed, for which *typical periods* were used as aggregation techniques [54, 76-79]. However, a common feature of these models is the existence of intraday storage technologies. As highlighted by Kotzur et al. [36], the implementation of storage technologies is only possible when the chronology within time steps is preserved, which is always the case for time steps within periods. For both, however, *typical time step* and *typical period* approaches have recently been developed to couple them to each other and enable storage modeling for *typical time steps* [21, 80] and of storage levels across multiple *typical periods* [68, 81].

This raises the question which of the two described aggregation approaches is superior and whether a general guideline can be provided as to how to aggregate input data to receive sufficient optimization results based on aggregated data.

# 3. Methodology

The following section introduces the comparison framework, aggregation methods, and models used for the validation of the developed temporal aggregation techniques, and introduces the data resources that were required to parametrize the models. In addition, the scientific purpose of each model's application with respect to an accurate temporal representation is discussed.

## 3.1. Comparison framework

The framework used to compare *typical periods* against *typical time steps* for ESOMs is shown in Figure 1. In order to enable a fair comparison to be made between the tested models and applied aggregation methods, we maintain the frameworks, aggregation intervals, and computing power constant for all models and model runs. We utilize the FINE[2] optimization framework [73] and time series aggregation module tsam[3] [10] to implement the models and aggregate the time series. All model runs were performed using the JURECA HPC Cluster [82] (for more details on the test system, see Section 3.3.4).

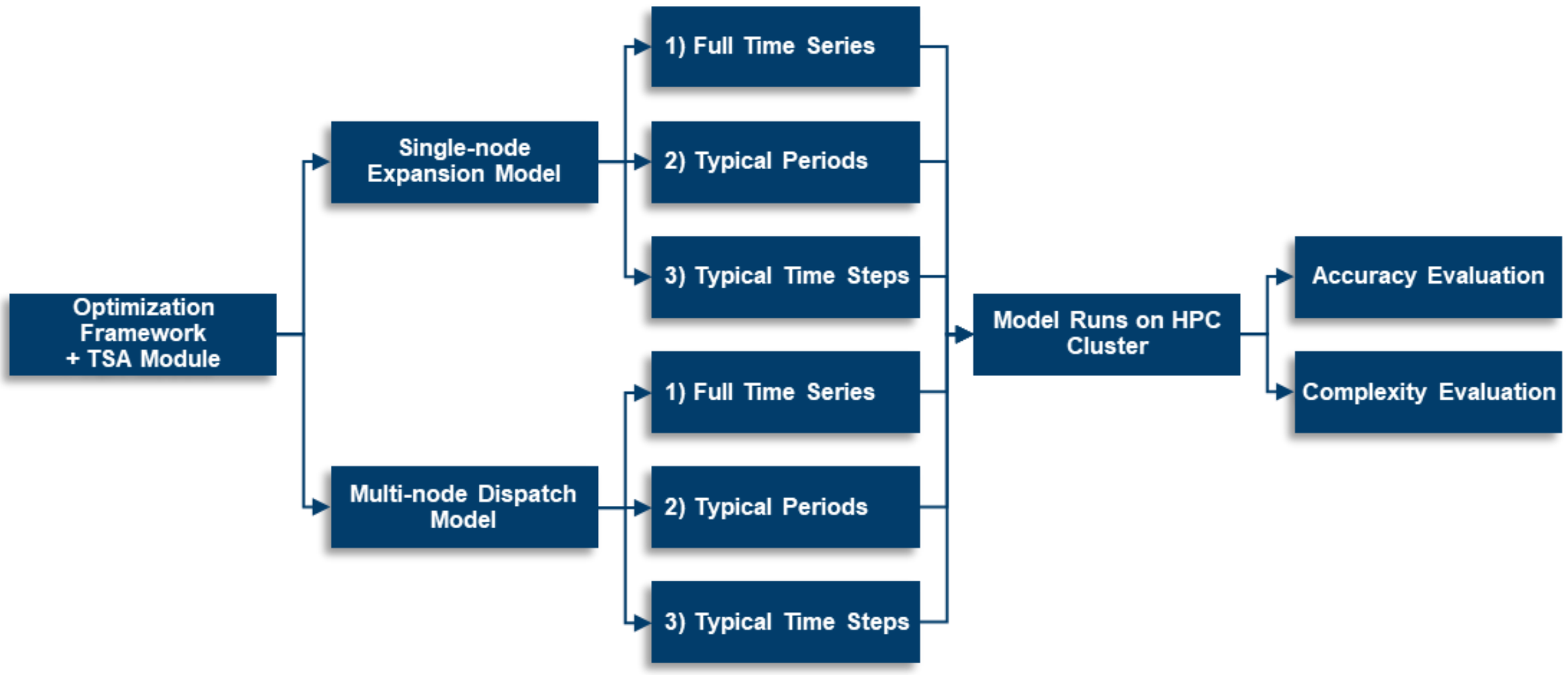

***Figure 1. Method applied to compare typical periods and typical time steps for different types of energy system optimization models.***

For all aggregation procedures, Ward's hierarchical clustering algorithm was chosen and each cluster of time steps was represented by its centroid. As the total number of considered time steps has a major impact on the computational complexity and to facilitate a comparison between *typical time steps* and *typical days*, only integer-multiples of 24 time steps were chosen for both *typical days* and *typical time steps*, which leads to the configuration of model runs and time steps per run shown in Table 2.

[2] FINE package available at: https://github.com/FZJ-IEK3-VSA/FINE (date accessed: 1/27/2021)

[3] tsam package available at: https://github.com/FZJ-IEK3-VSA/tsam (date accessed: 1/27/2021)

*Table 2. The model run configuration with respect to the total number of time steps considered.*

| Fully Resolved | Typical Days | Typical time steps |
|---|---|---|
| | | 24 x 1 |
| | | 48 x 1 |
| | | 96 x 1 |
| | 5 x 24 | 120 x 1 |
| | 10 x 24 | 240 x 1 |
| | 20 x 24 | 480 x 1 |
| | 40 x 24 | 960 x 1 |
| | 80 x 24 | 1920 x 1 |
| | 160 x 24 | 3840 x 1 |
| 8760 | 365 x 24 | 8760 x 1 |

This means that for *typical days*, 365 sample points comprising 24 dimensions of consecutive time steps per time series of the respective model were clustered, whereas for the *typical time steps*, 8760 data points with one dimension for each time series were aggregated. In the following section, the aggregation of the periods or time steps of multiple time series is explained in detail.

## 3.2. Temporal aggregation methods

The process of clustering-based temporal aggregation has been widely investigated in the literature. Although the fundamental steps to perform temporal aggregation of energy system models are the same, the algorithms utilized vary. Figure 2 summarizes the main steps that are included in most procedures to be found in the literature (c.f. [10, 19, 36, 83]).

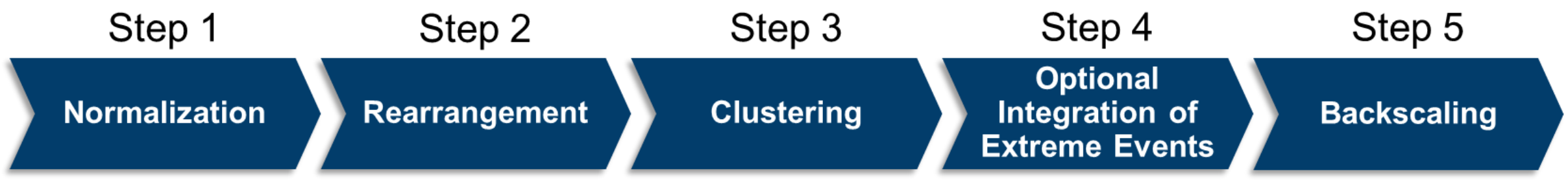

*Figure 2. Time series clustering procedure for energy system models.*

The input time series for bottom-up energy system optimization models generally consists of discrete time steps with a homogenous length, e.g., 1 h. In order to aggregate these time series by decreasing the number of time steps and likewise the mathematical complexity of the optimization problem, each time series $a \in \mathbb{A}$ comprises $s \in \mathbb{S}$ time steps.

As clustering generally incorporates all time series simultaneously, in the first step they are normalized to guarantee that no time series is overweighted in the clustering process due to its scale. In our case study, a min-max normalization was performed, i.e.:

$$x_{a,s} = \frac{x'_{a,s} - \min(x'_a)}{\max(x'_a) - \min(x'_a)}$$

In the second step, all normalized values of all time series $x_{a,s} \in \{0,1\}$ are rearranged into a set of candidate time steps in the case of *typical time steps* or candidate periods in the case of *typical periods*. In the case of *typical time steps*, the values are rearranged time step-wise, i.e., within a matrix that contains the time steps as rows and the attributes a columns:

$$\mathbf{D}_{ss} = (x_{a,s})^T = \begin{pmatrix} x_{1,1} & \cdots & x_{1,N_a} \\ \vdots & \ddots & \vdots \\ x_{N_s,1} & \cdots & x_{N_s,N_a} \end{pmatrix} = \begin{pmatrix} \mathbf{x}_{cand,1} \\ \vdots \\ \mathbf{x}_{cand,N_s} \end{pmatrix}$$

This can also be interpreted as a vector that contains the clustering candidates as samples. The rearrangement into candidate periods is more complex. For that, each time step within the time series is assigned a candidate time period $p \in \mathbb{P}$ and time step therein $t \in \mathbb{T}$, i.e.:

$$x_{a,s} = x_{a,p,t} \; with \; p = \left\lfloor \frac{s}{|\mathbb{T}|} \right\rfloor \; and \; t = s - p \times |\mathbb{T}|$$

For example, in the case of a time series with 8760 hourly time steps and a period length of 24 time steps, it would be 101. Time step would become the fifth of the fifth candidate period (i.e., the fifth hour of the fifth day), because:

$$p = \left\lfloor \frac{101}{24} \right\rfloor = 5 \; and \; t = 101 - 5 \times 24 = 5$$

This can be rearranged into the following candidate vector, in which each candidate is a row vector that contains $N_a \times N_t$ elements:

$$\mathbf{D}_{td} = \left(x_{a,p}\right)_t^T = \left(x_{p,a,t}\right) = \begin{pmatrix} x_{1,1,1} & \cdots & x_{1,1,N_t} & x_{1,2,1} & \cdots & x_{1,N_a,N_t} \\ \vdots & \ddots & \vdots & \vdots & \ddots & \vdots \\ x_{N_p,1,1} & \cdots & x_{N_p,1,N_t} & x_{N_p,2,1} & \cdots & x_{N_p,N_a,N_t} \end{pmatrix} = \begin{pmatrix} \mathbf{x}_{cand,1} \\ \vdots \\ \mathbf{x}_{cand,N_p} \end{pmatrix}$$

While the number of cluster candidates equals the number of time steps and the dimension of each candidate vector the number of attributes in the case of *typical time steps*, the number of cluster candidates equals the number of periods and the dimension of each candidate vector, the product of the number of attributes, and the number of time steps per period in the case of *typical periods*. Thus, for *typical time steps*, $|\mathbb{T}|$ times more candidates with a $|\mathbb{T}|$ times smaller dimension are clustered compared to *typical periods*.

In step 3, the respective cluster candidates are clustered according to a predefined number of clusters $k$ and a representative for each cluster $\mathbb{C}_k$ is calculated out of all candidates that are assigned to the respective cluster. Here, the cluster algorithm used and the method for determining a representative may differ. In this study, Ward's hierarchical clustering algorithm was used and each cluster was represented by its centroid, i.e., in the case of *typical time steps*:

$$c_{k,a} = \frac{1}{|\mathbb{C}_k|} \sum_{s \in \mathbb{C}_k} x_{s,a}$$

And in the case of *typical periods*:

$$c_{k,a,t} = \frac{1}{|\mathbb{C}_k|} \sum_{p \in \mathbb{C}_k} x_{p,a,t}$$

Step 4 is optional, and one in which the extreme values of the original time series can be manually added so that the cluster representatives are rescaled to satisfy certain features of the original time series, e.g., its mean value. In order to not distort the comparison between *typical time steps* and *typical days*, this step is skipped herein this work. Furthermore, the mean values of the original time series are retained because of the clusters' representations by centroids.

In step 5, the clustered time series are finally rescaled, i.e., in the case of *typical time steps*:

$$c'_{k,a} = \min(x'_a) + c_{k,a}(\max(x'_a) - \min(x'_a))$$

And in the case of *typical periods*:

$$c'_{k,a,t} = \min(x'_a) + c_{k,a,t}(\max(x'_a) - \min(x'_a))$$

Finally, but importantly, each time-dependent cost related to an aggregated time step or period must be weighted in the aggregated model with respect to its number of occurences $|\mathbb{C}_k|$. Furthermore, the *typical time steps* and *typical periods* are coupled in the subsequent model runs in order to model energy storage, as presented by Wogrin et al. [21] for *typical time steps* and Kotzur et al. [81] for seasonal storage based on *typical periods*.

## 3.3. Case study

For our analysis, we employ two fundamentally different models, namely a single-node model for expansion planning and a multi-node dispatch model. The single-node model optimizes the choice, scale, and operation of energy system technologies for a self-sufficient building including multiple commodities and storage technologies; the multi-node model optimizes the dispatch schedule for the German electricity system.

The models are meant to represent extreme cases of potential modeling applications: In the case of the building model, the isolation of the energy system leads to a significant cost contribution of storage technologies and requires the temporal aggregation technique to remain flexible with respect to different storage cycles in order to enable realistic storage sizing. On the other hand, the dispatch model does not explicitly consider storage, but hundreds of multiregional time series, which could result in fundamentally different requirements for the aggregation technique.

Both models have been outlined in earlier works by the authors, but are briefly introduced below to provide an overview of their respective features. For a more detailed description, see Kotzur et al. [84] for the single-node model and Priesman et al. [6] for the multi-node one.

### 3.3.1. The Self-Sufficient Building Model

The first case study, a model of a self-sufficient building, is an island system that exclusively relies on a RES, namely solar photovoltaic (PV). It must not only fit electricity demand but also heat demand. This leads to a complex setup and requires higher storage capacities, as the system lacks a freely dispatchable energy source. Therefore, a wide set of potential technologies is considered in this model, including electric, thermal, hydrogen, and liquid organic hydrogen carrier (LOHC) storage units. A simplified scheme of the model based on a layout by Kotzur et al. [84] that was transferred into the FINE framework is displayed in Figure 3.

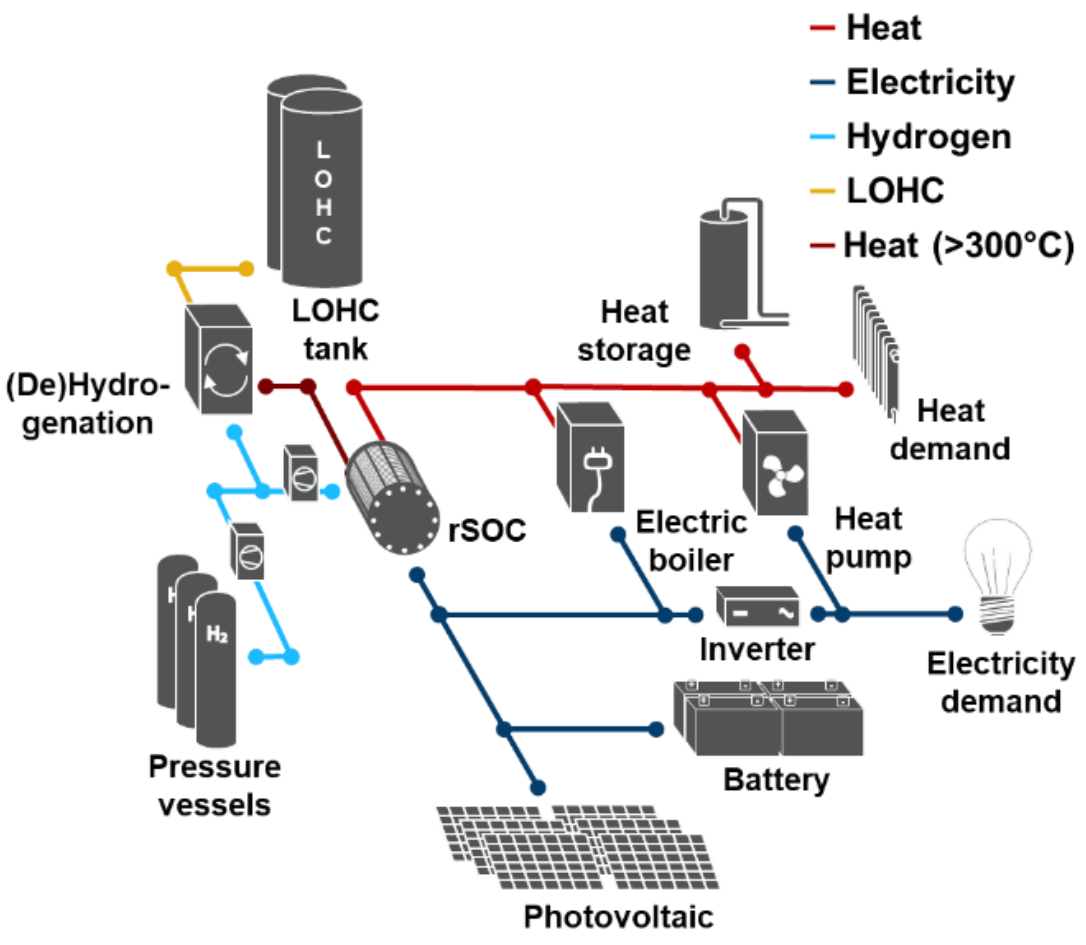

*Figure 3. A simplified scheme of a self-sufficient residential energy system.*

The self-sufficient building model provides numerous options for meeting the building's heat and electricity demand. Furthermore, it is clear that the energy storage components will most likely work on different time scales, depending on the form of energy that they are storing.

Table 3 provides an overview of the components' cost contributions and the data is qualitatively discussed with respect to its meaning for the model from a mathematical point of view. The exact manifestation of the assumptions is not relevant to the research question, and is therefore not discussed.

***Table 3. Cost parameters of the self-sufficient building model.***

| | Capex | | | | Opex | | Lifetime | | Source |
|---|---|---|---|---|---|---|---|---|---|
| Components | Fixed | | Variable | | Fixed+Variable | | | | |
| Photovoltaic Ground | — | — | 4000 | €/kW$_{p}$ | 1 | % Inv./a | 20 | a | |
| Photovoltaic Rooftop | — | — | 769 | €/kW$_{p}$ | 1 | % Inv./a | 20 | a | [85] |
| Inverter | — | — | 75 | €/kW$_{p}$ | — | — | 20 | a | [86] |
| Battery | — | — | 301 | €/kWh$_{p}$ | — | — | 15 | a | [85] |
| Reversible Solid Oxide Cell | 5,000 | € | 2,400 | €/kW$_{el}$ | 1 | % Inv./a | 15 | a | [87] |
| Heatpump | 4,230 | € | 504.9 | €/kW$_{th}$ | 1.5 | % Inv./a | 20 | a | [88] |
| Thermal Storage | — | — | 90 | €/kWh$_{th}$ | 0.01 | % Inv./a | 25 | a | [89] |
| E-Heater & E-Boiler | — | — | 60 | €/kW$_{th}$ | 2 | % Inv./a | 30 | a | [89] |
| Tank | — | — | 0.79 | €/kWh$_{H2}$ | — | — | 25 | a | [90] |
| Dibenzyltoluene | — | — | 1.25 | €/kWh$_{H2}$ | — | — | 25 | a | [91, 92] |
| Hydrogen Vessels | — | — | 15 | €/kWh$_{H2}$ | — | — | 25 | a | [93] |
| Hydrogenizer | 2,123.3 | € | 761.1 | €/kW$_{H2}$ | 1 | % Inv./a | 20 | a | [91] |
| Dehydrogenizer | 1,140 | € | 408.6 | €/kW$_{H2}$ | 1 | % Inv./a | 20 | a | [91] |
| Low Pressure Compressor | — | — | 1716.71 | €/kW$_{p}$ | 1 | % Inv./a | 25 | a | [94] |
| High Pressure Compressor | 560 | € | 1329.8 | €/kW$_{p}$ | 1 | % Inv./a | 25 | a | [94] |
| Heat-Exchangers 1 and 2 | — | — | 1 | €/kW$_{th}$ | 1 | % Inv./a | — | a | — |
| Expanders 1 and 2 | — | — | 1 | €/kW$_{th}$ | 1 | % Inv./a | 25 | a | — |

The time series data of the building model comprises five different time series with hourly resolutions: Three solar profiles for the three differently-oriented PV modules, as well as the electricity and heat demand time series. The profiles were generated using the open source software tsib [95, 96], which was used to transform weather information, such as global horizontal irradiance and ambient temperature, into demand profiles for archetype buildings and capacity factors for renewable energy sources used in the domestic sector. The raw weather data was based on the COSMO dataset [97-99] pertaining to the city of Berlin for the year 2013.

### 3.3.2. The Electricity Dispatch Model

The Electricity Dispatch Model does not consider capacity expansion, i.e., the size of all components is predefined and only their operation is optimized. The multi-regional model covers 16 regions based on the NUTS1 regions of Germany, which correspond to the country's federal states. Accordingly, power flows between the different regions, as well as neighboring countries, are considered. This as-pect of the model is an interesting feature with respect to the impact of temporal aggregation on spatially-resolved data, as the aggregation of temporal data not only affects the behavior of single components but also the interaction of components at different locations, as demonstrated by Caglayan et al. [100].

The electricity dispatch model is a simplified FINE implementation of the JERICHO dispatch model developed and published by Priesmann et al. [6]. The model is schematically illustrated in Figure 4.

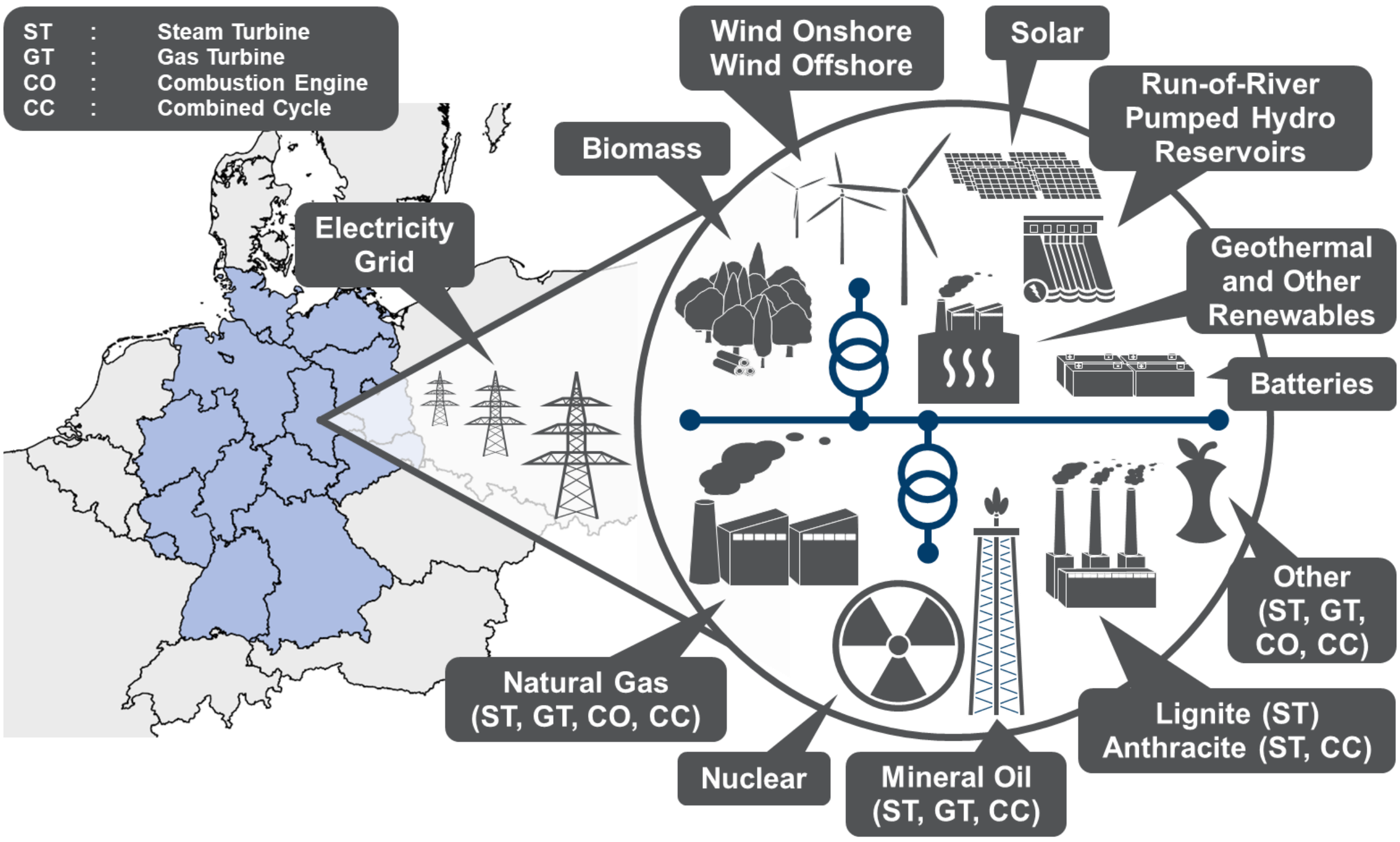

*Figure 4. Layout of the dispatch model.*

The dispatch model incorporates 25 different electricity-generating technologies, of which 15 are based on fossil fuels, eight on RES and two are storage technologies. For the analyses conducted in this study, the storage dispatch was fixed prior to the optimization on the basis of historical values. The technologies based on fossil fuels, as well as biomass, i.e., the freely-dispatchable technologies, were modeled as conversion units that consume a certain amount of lignite, hard coal (anthracite), uranium, methane, or biomass to produce an amount of electricity and $CO_2$ emissions. The non-dispatchable technologies, comprising wind and solar energy and run-of-the-river hydroelectricity, are modeled as sources of fixed capacity but with constraining capacity factors at an hourly resolution.

The objective of this optimization model is to minimize the costs of covering the inelastic electricity demand at each hour. The costs are the sum of the marginal costs for operating the conversion technologies and based on the time-dependent prices for fuels and EU allowances (EUAs) within the EU Emissions Trading Scheme and the technology-specific operation costs.

The operation of all units is strictly linearly modeled, i.e., unit commitment constraints such as minimum up- and down-time or startup and shutdown costs that necessitate binary variables are not considered. This leads to the generation of a large, but comparably easily solvable energy system model.

As already noted, the dispatch model considers neither capacity expansion, nor depreciation costs, which is the most significant difference to the self-sufficient building model. On the other hand, it comprises multiple regions and thus transmission grids. In order to provide an overview of the model's technological structure, Figure 5 depicts the share of each electricity supply technology at each of the 16 NUTS1 regions in Germany, as well as the inner German transmission lines and those to neighboring countries, offering the option to export electricity to Germany. Furthermore, the cumulative yearly electricity consumption in each region is highlighted by a corresponding background color.

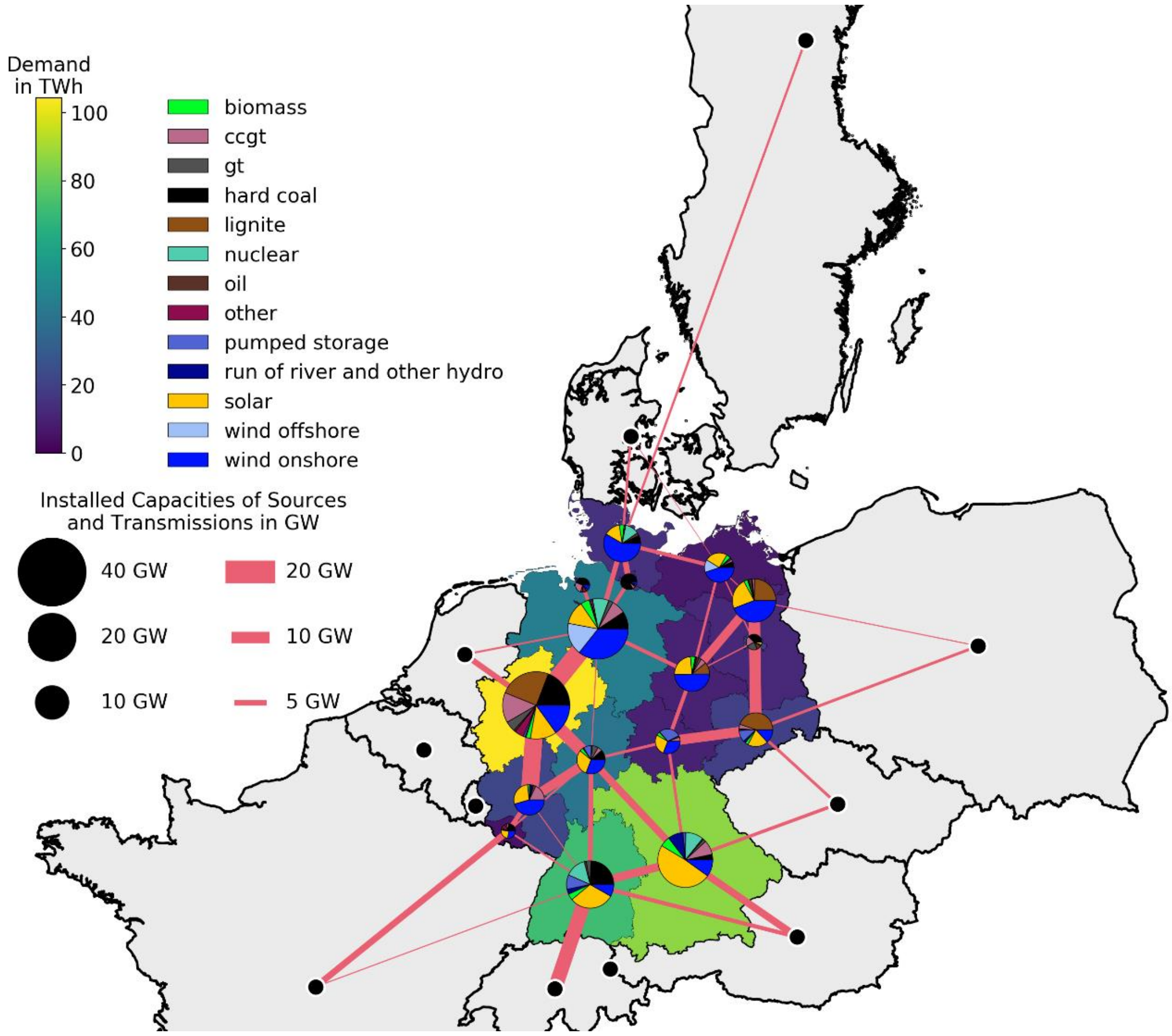

*Figure 5. Input data of the dispatch model: Aggregated transmission lines, installed capacities, and yearly energy demand*

The list of plant technologies, their rated power, and position are drawn from the Federal Network Agency (Bundesnetzagentur) for 2019 [101] as well as the Internation Energy Agency (IEA) and Nuclear Energy Agency (NEA) [102]. The data for the length and capacity of the inner- and trans-German transmission lines are derived from the preceding project SciGRID [103]. The geoinformation data used for the illustrations was obtained from Eurostat [104]. The electricity demand time series are taken from the ENTSO-E Transparency Platform [105] and regionalized using data from the Country Working Group on Energy Balances (Länderarbeitskreis Energiebilanzen) [106] and the Federal and State Statistical Offices (Statistische Ämter des Bundes und der Länder) [107]. Finally, the import costs per country and the capacity factors for renewable energy sources were all derived from the ENTSO-E database for the year 2019 [105]. The techno-economic assumptions of the generation units are depicted in Table 4.

*Table 4. Cost parameters of the dispatch model.*

| Components | Efficiency | | Operational cost | | Emission factor | | Number of units |
|---|---|---|---|---|---|---|---|
| **Biomass** | 38-45 | % | 22.2-23.0 | €/MWh | 0.016 | $t_{CO_2}/MWh$ | 16,157 |
| **Gas (combined cycle)** | 41-63 | % | 5.0-8.7 | €/MWh | 0.204 | $t_{CO_2}/MWh$ | 90 |
| **Gas (simple cycle)** | 29-42 | % | 4.2-7.0 | €/MWh | 0.204 | $t_{CO_2}/MWh$ | 184 |
| **Hard coal** | 23-50 | % | 10.6-16.8 | €/MWh | 0.342 | $t_{CO_2}/MWh$ | 77 |
| **Lignite** | 20-43 | % | 11.1-19.3 | €/MWh | 0.4 | $t_{CO_2}/MWh$ | 60 |
| **Nuclear** | 33 | % | 7.5 | €/MWh | 0 | $t_{CO_2}/MWh$ | 7 |
| **Mineral Oil** | 30-40 | % | 4.8-8.4 | €/MWh | 0.266 | $t_{CO_2}/MWh$ | 47 |
| **Hydro (pumped storage)** | 75 | % | 10.0 | €/MWh | — | — | 74 |
| **Hydro (run of river)** | — | — | 10.0 | €/MWh | — | — | 7,254 |
| **Photovoltaic (ground)** | — | — | 7.5 | €/MWh | — | — | 8,795 |
| **Photovoltaic (rooftop)** | — | — | 7.5 | €/MWh | — | — | 1,991,267 |
| **Wind (onshore)** | — | — | 7.5 | €/MWh | — | — | 28,777 |
| **Wind (offshore)** | — | — | 17.0 | €/MWh | — | — | 1,497 |
| **Other (generation)** | 15-45 | % | 6.1-7.0 | €/MWh | 0.016 – 0.3 | $t_{CO_2}/MWh$ | 463 |
| **Other (storage)** | 90 | % | 5.0-10.0 | €/MWh | — | — | 112,154 |

### 3.3.3. Accuracy indicators

In order to quantify the quality of the cluster-based temporal aggregation for energy system models, different metrics have been established in the literature, among which, the Root Mean Squared Error (RMSE) between the original and aggregated time series [19, 36] and the RMSE of the respective duration curves [13, 19, 40, 42, 72, 83, 108-111] have gained popularity. Furthermore, the MAE traditionally provides an alternative to the RMSE, and is therefore also considered.

In the following, we define the RMSE of an aggregated time series as follows:

$$\mathrm{RMSE} = \sqrt{\frac{1}{|\mathbb{P}| \times |\mathbb{T}|}\sum_{|\mathbb{C}|}\sum_{p\in\mathbb{C}_k}\sum_{|\mathbb{T}|}\left(x_{p,t}-\tilde{x}_{k,t}\right)^2}$$

where $|\mathbb{P}|$ is the number of periods, $|\mathbb{T}|$ is the number of time steps per periods, and $|\mathbb{C}|$ is the number of representative periods (*typical days* or *typical time steps*). $x_{p,t}$ is the value of the original time series at a specific time step, e.g., for *typical days* with an hourly resolution, $x_{3,6}$ is the value of a time series on January 3$^{rd}$ at 6 a.m., and $\tilde{x}_{k,t}$ is the value of the $k^{th}$ cluster, which represents the respective value $x_{p,t}$ in the aggregated time series. For multiple time series (attributes) $a$, we form the RMSE as follows:

$$\mathrm{RMSE}_{\mathrm{tot}} = \sqrt{\frac{1}{|\mathbb{A}| \times |\mathbb{P}| \times |\mathbb{T}|}\sum_{|\mathbb{A}|}\sum_{|\mathbb{C}|}\sum_{p\in\mathbb{C}_k}\sum_{|\mathbb{T}|}\left(x_{a,p,t}-\tilde{x}_{a,k,t}\right)^2} = \sqrt{\frac{1}{|\mathbb{A}|}\sum_{|\mathbb{A}|}\mathrm{RMSE}_a^2}$$

In order to calculate the RMSE of the aggregated time series' duration curves, the values of the clusters' representatives are repeated according to the respective cluster size $|\mathbb{C}_k|$, i.e., the number of days that a typical day is representive of, or the number of time steps a system state represents. If the values of the original and predicted time series are sorted, we obtain the duration curves $x_{DC,s}$ and $\tilde{x}_{DC,s}$. Accordingly, the RMSE of a duration curve is calculated as follows:

$$\mathrm{RMSE_{DC}} = \sqrt{\frac{1}{|\mathbb{S}|}\sum_{|\mathbb{S}|}\left(x_{DC,s} - \tilde{x}_{DC,s}\right)^2}$$

where $|\mathbb{S}|$ is the total number of time steps. For multiple time series, the RMSE of the duration curves is analogously calculated:

$$\mathrm{RMSE_{DC,tot}} = \sqrt{\frac{1}{|\mathbb{A}| \times |\mathbb{S}|}\sum_{|\mathbb{A}|}\sum_{|\mathbb{S}|}\left(x_{DC,a,s} - \tilde{x}_{DC,a,s}\right)^2} = \sqrt{\frac{1}{|\mathbb{A}|}\sum_{|\mathbb{A}|}\mathrm{RMSE}^2_{DC,a}}$$

The MAE is the final accuracy indicator, and uses the average error instead of the squared error, i.e.:

$$\mathrm{MAE} = \frac{1}{|\mathbb{P}| \times |\mathbb{T}|}\sum_{|\mathbb{C}|}\sum_{p\in\mathbb{C}_k}\sum_{|\mathbb{T}|}\left|x_{p,t} - \tilde{x}_{k,t}\right|$$

and

$$\mathrm{MAE_{tot}} = \frac{1}{|\mathbb{A}| \times |\mathbb{P}| \times |\mathbb{T}|}\sum_{|\mathbb{A}|}\sum_{|\mathbb{C}|}\sum_{p\in\mathbb{C}_k}\sum_{|\mathbb{T}|}\left|x_{a,p,t} - \tilde{x}_{a,k,t}\right|$$

### 3.3.4. Computational test system

The models were run on a computing node with the configurations depicted in Table 5. Each run was carried out on a core with one thread for workload orchestration and one for performing the optimization. Once a core had completed a taks, it was assigned to a new model optimization. To address our research question, we performed a cross-comparison of both aggregation methods, i.e., using *typical days* and *typical time steps*. To this end, we conducted several runs with both models and benchmarked the results to those of the fully resolved cases.

***Table 5. Specifications of the computing node used for the optimization runs.***

| | |
|---|---|
| CPU model | Intel(R) Xeon(R) CPU E5-2697 v3 |
| Number of cores per computing node | 28 |
| Threads per core | 2 |
| CPU max frequency [MHz] | 3,600 |
| Shared Memory [GB] | 1,024 |

# 4. Results and Discussion

The following section first investigates the impact of aggregated configurations with respect to the accuracy of the data aggregation using a set of clustering indicators. Subsequently, we compare the deviation of the aggregated models' optimal objective to the fully resolved reference cases, with consideration to the individual computational speedup. Then, the aggregation-induced deviations of both models are analyzed in detail by investigating the cost contributions of the systems' components. Finally, the time consumption of the individual processes, as well as the memory consumption of each model, is evaluated over time.

## 4.1. Input data-driven analysis of the aggregation accuracy

Figure 6 shows the considered accuracy indicators (see Section 3.3.3) for the normed time series (with values between 0 and 1) of the self-sufficient building (top) and dispatch model (bottom) for both *typical periods* and *typical time steps*.

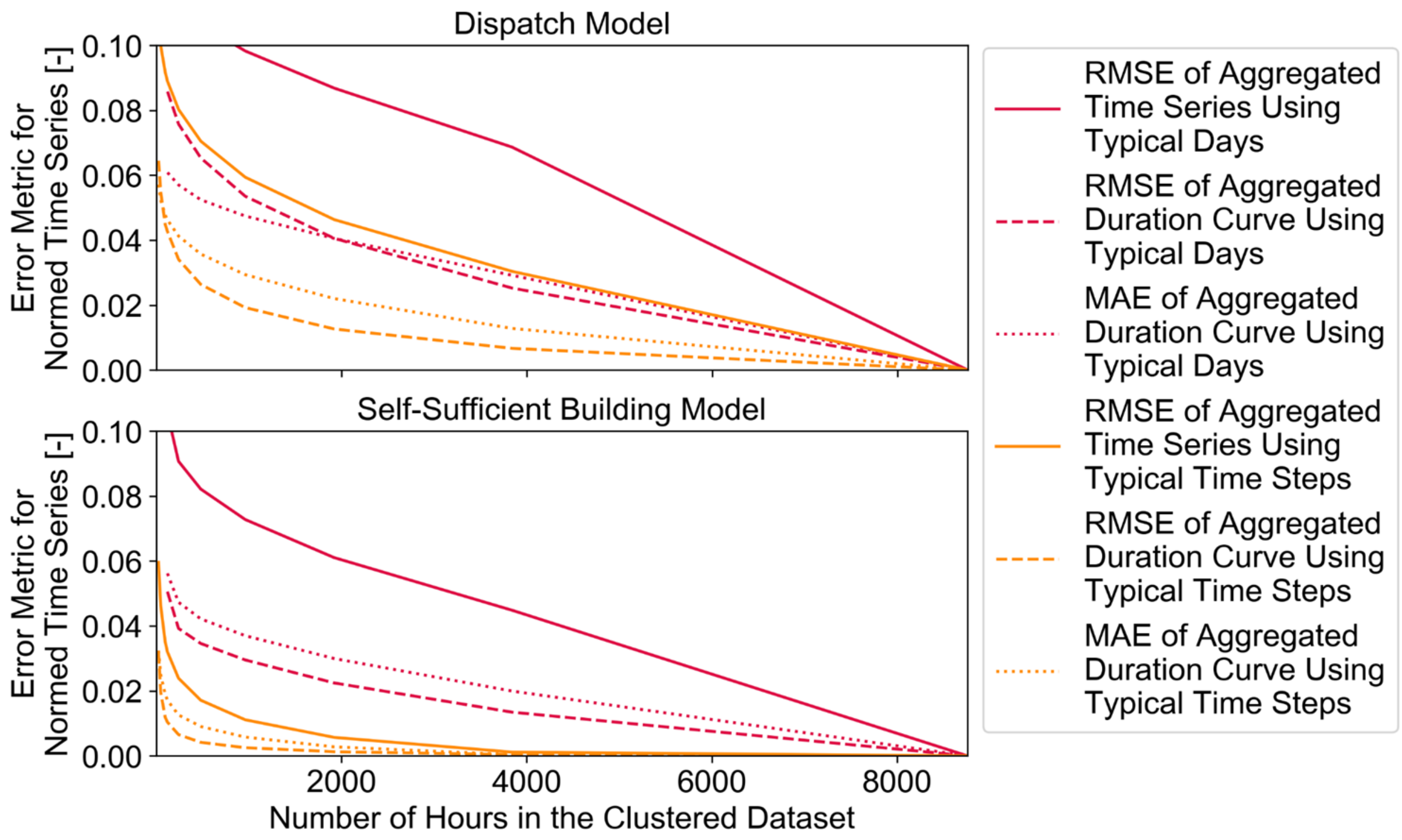

***Figure 6. Accuracy indicators for the normed time series (values of between 0 and 1) of the self-sufficient building (top) and dispatch model (bottom) for both typical days and typical time steps.***

The x-axis shows the number of equivalent *typical time steps*, i.e., the number of time steps to which the time series has been aggregated. In the case of *typical days*, this equals the number of *typical days* multiplied by 24 time steps per day. Here, it can be seen that, from an input data perspective, the aggregation of the time series to single time steps outperforms the aggregation to conjoined periods for any clustering indicator and both models. The reason for this is that in the case of *typical time steps*, each time step is a single candidate for clustering and can be freely grouped with any other time step, which is similar enough. When *typical days* are clustered, however, time steps are only compared to other time steps in the same daytime. Furthermore, the clustering of *typical days* is performed in a 24-times higher dimensional space with 24-times fewer candidates. Therefore, the 'curse of dimensionality' leads to higher aggregation-induced deviations from the original time series.

## 4.2. Output data-driven analysis of the aggregation accuracy

The main purpose of the aggregation techniques for energy system optimization models is to reduce computational complexity while keeping the impact on the optimized solution as small as possible. As the objective function of the considered models comprises the total annual costs (TAC), the main determinants for the quality of a good aggregation are the calculation time and a small deviation from the original, optimal objective.

Figure 7 shows the absolute deviation from the optimized objective function value from the reference case over the calculation time for the self-sufficient building (left) and dispatch model (right).

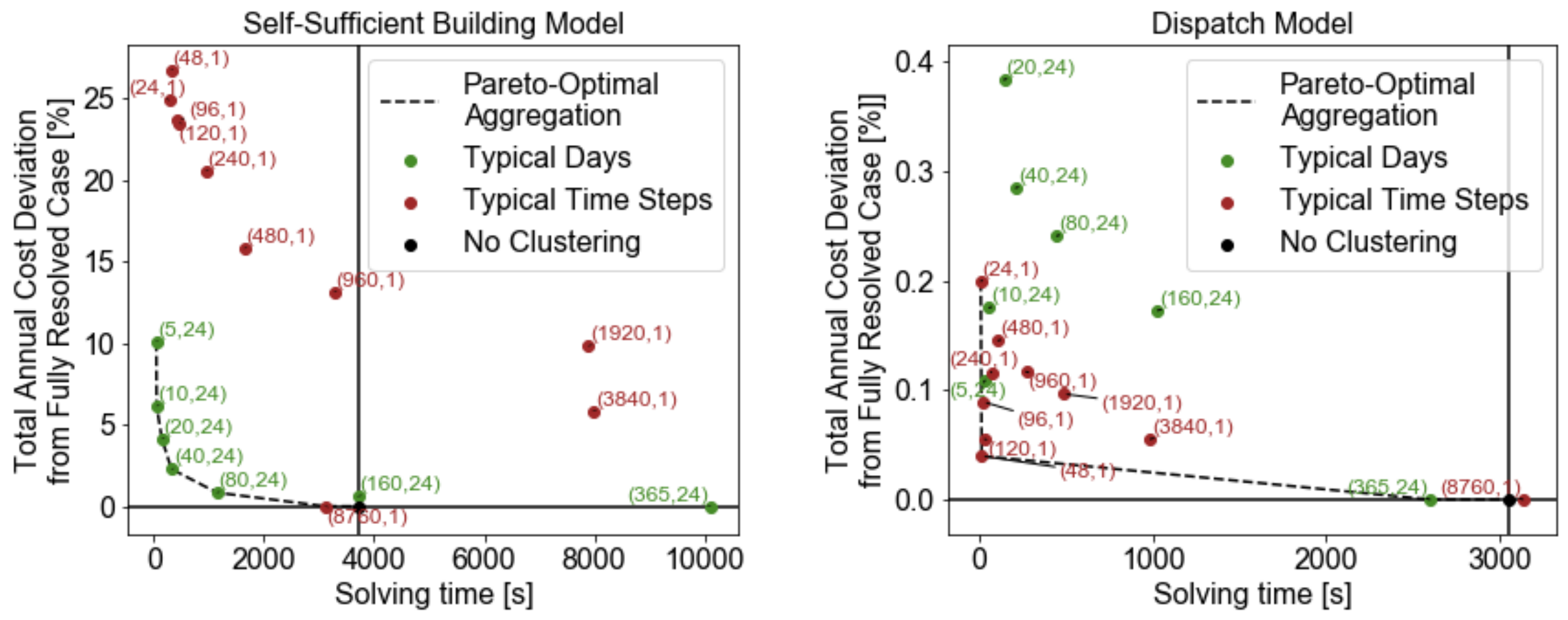

***Figure 7. Absolute deviation from the optimized objective function value of the fully resolved case over the calculation time for the self-sufficient building (left) and dispatch model (right).***

In the case of the self-sufficient building model, the deviation from the optimal objective is consistently much greater than in the case of the dispatch model. Accordingly, it is much more sensitive to temporal aggregation due to the strong temporal coupling of the fully-resolved model, as induced by the numerous storage technologies. Although *typical time steps* consistently outperform the aggregation to *typical days* with respect to the clustering indicators, and therefore the deviation between the aggregated and original time series, it is evident from Figure 7 that the *typical time steps* represent the superior aggregation option for the dispatch model, as they lead to smaller deviations with low calculation times. In contrast, they also lead to significantly larger deviations from the optimal objective value of the fully resolved self-sufficient building model. Here, *typical days* are the better aggregation option, which is also supported by the dashed lines in Figure 7 representing the respective pareto fronts. Furthermore, it can be observed that the optimal objectives converge to the one of the fully resolved case if the number of *typical time steps* or *typical days* is increased, which can be expected, but which supports the mathematical validity of the aggregation. Moreover, it can be observed that, especially in the case of the self-sufficient building, the modeling approach for seasonal storage with a reduced number of capacity constraints introduced by Kotzur et al. [81] leads to a computational overhead in the case of aggregation configurations with a large number of *typical days* or *typical time steps*. Similar observations were made by Wogrin et al. [21] with respect to temporally-linked *typical time steps* to model storage components, in spite of the aggregation.

In order to further analyze the aggregation-induced errors, the individual cost contributions for the self-sufficient model (top) and dispatch model (bottom) for all aggregation configurations are depicted in Figure 8.

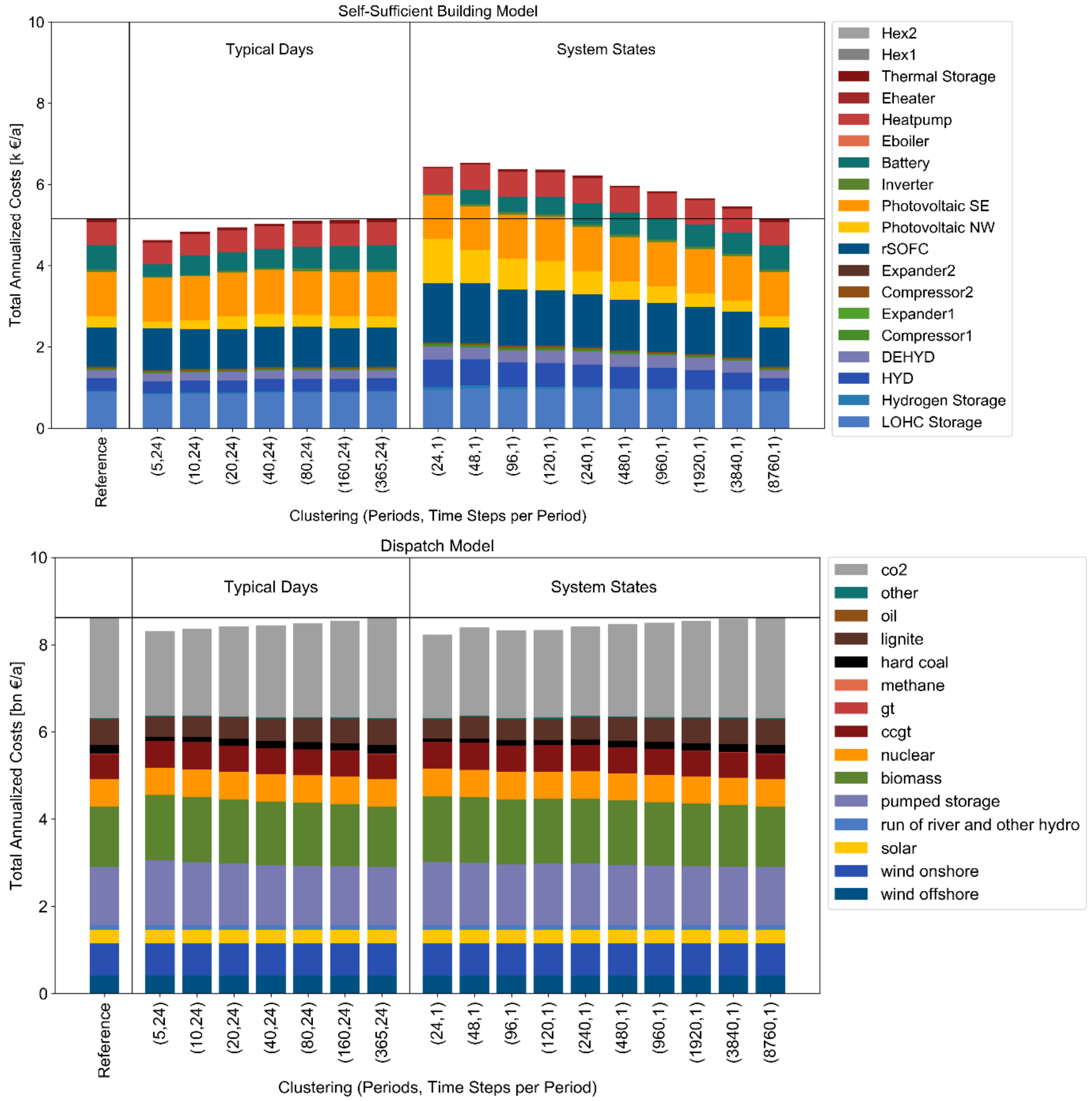

***Figure 8. Individual cost contributions for the self-sufficient model (top) and dispatch model (bottom) for all aggregation configurations.***

As the Figure 8 shows, the centroid-based aggregation to *typical time steps* or *typical days* generally leads, in three out of four cases, to a consistent underestimation of the reference optimal objective values. In the case of *typical time steps* for the self-sufficient building, however, the aggregation leads to a significant overestimation of the optimal objective value. Interestingly, this observation approves the inconsistent results in the literature, e.g., Zatti et al. [71] observed that, depending on the model, the centroid-based aggregation of the k-means algorithm could lead to over- or underestimation. This indicates that the impact of the aggregation method is highly dependent on the model to which it is applied.

In general, the underestimation of the majority of the aggregated models can be attributed to the fact that the time series are affected by the averaging effect of the centroid-based cluster representation, i.e., the minimum values of a time series are overestimated by the aggregation, whereas the maximum values are underrestimated. As is mathematically shown by Teichgräber et al. [112], the centroid-based representation of the time series that form coefficients in the constraint vector lead to a relaxation of temporally-decoupled models. This applies, for instance, to demand time series. However, in the case of the dispatch model, the cost time series, which yield coefficients in the cost vector of the optimization problem, are clustered as well. As Teichgräber et al. [112] have shown, the aggregation of the cost vector and the reduction in the corresponding variables reduce the optimization model's feasible solution space. However, this counteracting effect is of minor importance for the dispatch model, as the objective function based on *typical days* or *typical time steps* is consistently smaller than the one of the reference case. However, the impact of the latter effect can be observed between 48 and 96 *typical time steps* in the case of the dispatch model shown in Figure 8. Here, the total annual costs decrease for a higher temporal resolution, despite the use of hierarchical clustering, which should gradually decrease the aggregation-induced error. This implies the existence of at least two opposing effects on the optimal objective value due to aggregation: One that increases the optimal objective due to fewer variables in the aggregated model and one that decreases it due to fewer constraints.

Storage components, as in the case of the self-sufficient building, lead to a much less predictable impact of centroid-based clustering compared to the dispatch model, as shown in the upper half of Figure 8. Compared to *typical time steps*, the aggregation using *typical days* has a relatively small impact on the component's cost contributions. Yet, it can be observed that a smaller number of typical days most prominently decreases the optimum battery capacity, as the aggregation-induced smoothing of the aggregated time series profiles intensifies, and therefore smaller battery capacities are needed to balance intra-daily demand and supply fluctuations.

In contrast to that, an increase in the number of *typical time steps* decreases the capacity-induced cost contribution of the northwestern-oriented PV panels and the reversible solid oxide cell (rSOC). Simultaneously, the size of the battery slowly increases with a higher number of *typical time steps*.

The battery has relatively large capacity-specific investment costs (€/$kWh_{storage}$) and a relatively high self-discharge rate, whereas the hydrogen storage has higher power-specific investment costs (€/$kW_{(dis)charge}$) due to the expensive rSOC. This predetermines batteries for short daily storage cycles and small capacities, whereas hydrogen storage is preferable for large capacities and long storage cycles at low charge or discharge rates.

In the case of *typical time steps*, all storage components function identically at time steps that are assigned to the same cluster based on their similarities. The fact that the electricity demand profile and solar profiles resemble themselves during the morning and evening hours leads to an eventual assignment of these hours to the same cluster. This in turn forces the battery to operate in the same way during morning and evening hours and inhibits a daily operation cycle. Accordingly, in an aggregated model, it is less economically-viable based on *typical time steps*. Instead, the operational feasibility of the model is maintained by oversizing the solar panels and hydrogen subsystem used for storing energy over longer time scales.

In summary, the existence of storages is one of the most determinant factors for an accurate temporal aggregation, even if the time series have a strong daily pattern (which is also the case for the dispatch model).

### 4.3. In-depth analysis of computational resources

Apart from the manifold impacts that temporal aggregation has on a model's optimal solution, the main purpose of temporal aggregation remains computational speedup. Therefore, the contributions of the different processes to the overall runtime, as well as individual memory consumption, are analyzed in the following. Figure 9 illustrates the runtime contributions of each process for the self-sufficient building model (top) and dispatch model (bottom) for all aggregation configurations.

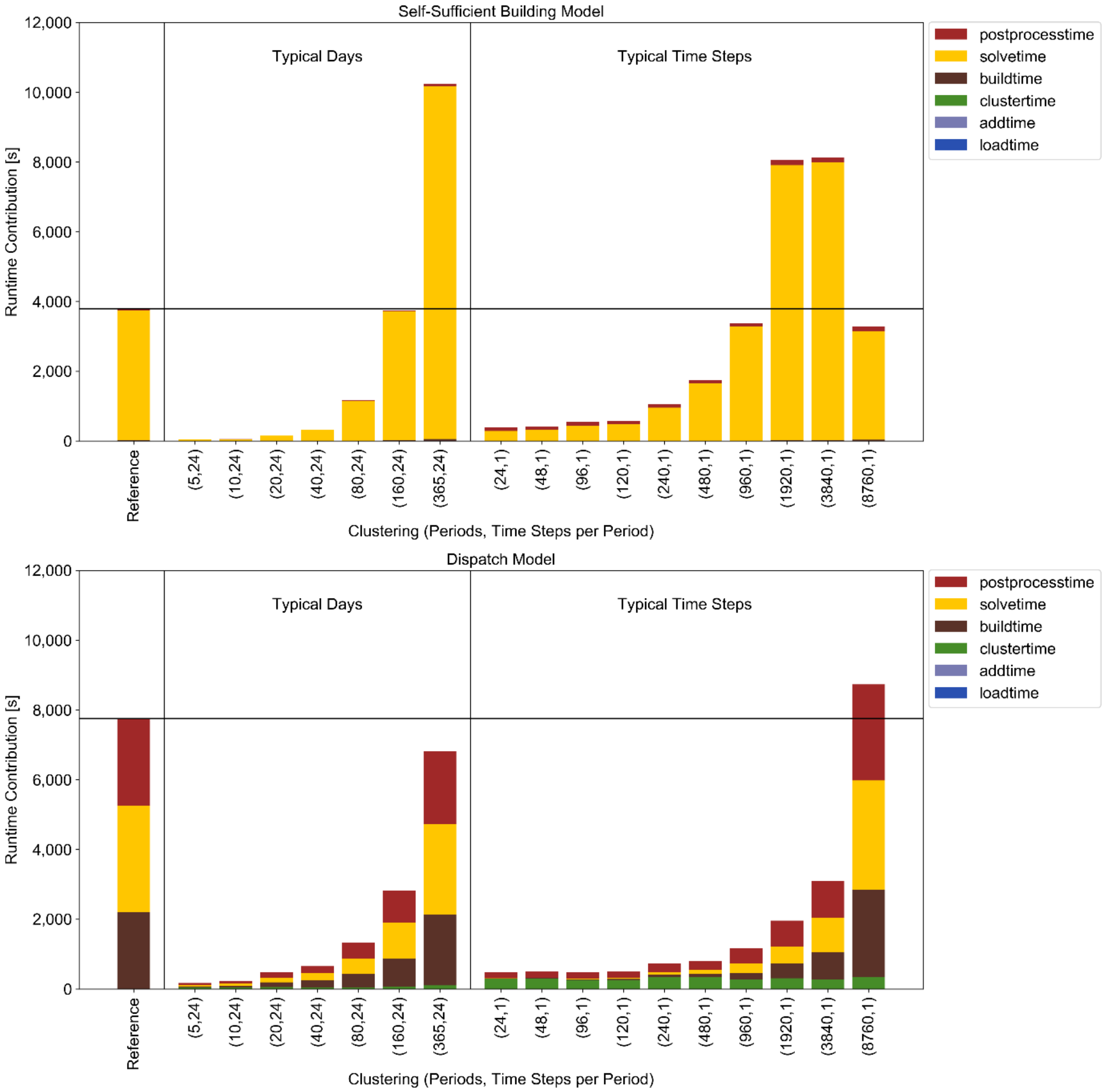

***Figure 9. Calculation times of the self-sufficient building (top) and dispatch model (bottom), depending on the aggregation configuration.***

As can be seen in the case of the self-sufficient building, almost all of the overall runtime is needed to solve the optimization problem. The time consumption for loading the data, adding the components, clustering the input data, setting up the mathematical problem structure using Pyomo, and mapping the solution from Gurobi back to the model's variables, is negligible. This can be explained by the mathematical connectivity [7] of the variables and constraints caused by a complex component structure and time step linking storage equations whereas, on the other hand, only one region and five time series are considered.

In contrast, the overall runtime of the dispatch model is primarily determined by three operations: The time to set up the problem structure, the time to solve the problem, and the time to map the solution back to the set-up problem's structure. The reason for this is the large amount of input data given by 6325 time series, as well as a large number of regions (16 national and nine international), whereasthe problem structure is relatively simple because of the temporal decoupling due to the absence of storage components and a flat system structure, because all of the components are connected to a single electricity grid. Moreover, the time for clustering becomes a significant share of the overall calculation time in the case of aggregating *typical time steps*, which can be attributed to the fact that the computational complexity of Ward's hierarchical clustering algorithm scales with $O(n^2)$, and that *typical time steps* feature 24 times more starting candidates than *typical days*. Moreover, the process itself is much slower than for the self-sufficient building due to the large amount of time series. Apart from that, the aggregation tends to take longer for stronger aggregated model configuration, which makes the aggregation time an even more significant runtime contribution for strongly aggregated models. On the other hand, this motivates the use of simple clustering approaches instead of new but computationally-expensive aggregation techniques that erode the original purpose of temporal aggregation, i.e., the computational speedup.

Figure 10 shows the memory consumption of the self-sufficient building (top) and dispatch model (bottom), depending on the aggregation configuration.

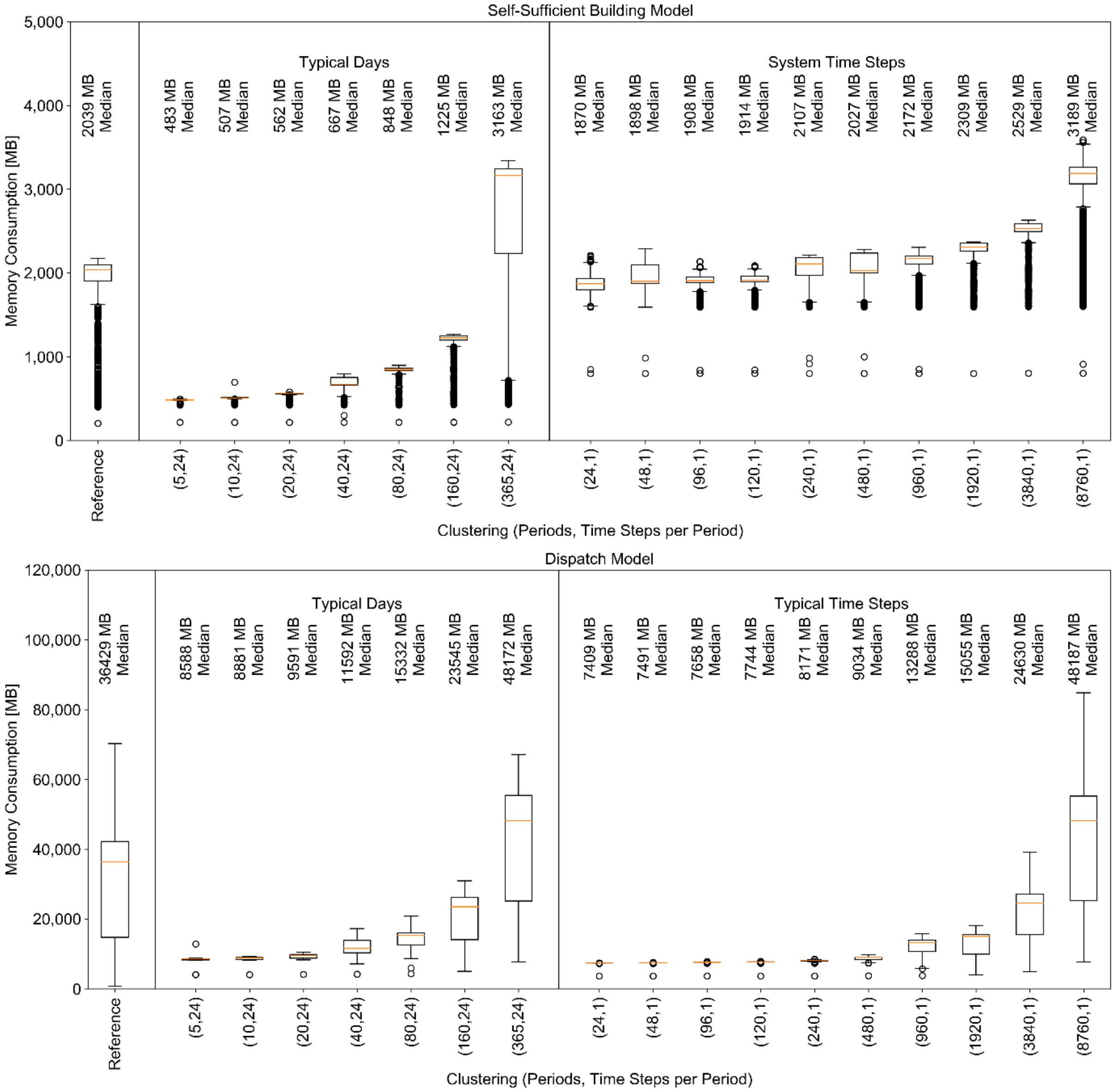

***Figure 10. Memory consumption of the self-sufficient building (top) and dispatch model (bottom) depending on the aggregation configuration.***

Another interesting difference between both models is shown in Figure 10: As expected, the memory consumption of both models increases with a higher number of *typical days* or *typical time steps*. However, the total amount of memory consumption of the dispatch model is more than an order of magnitude larger than that of the corresponding aggregation configuration of the self-sufficient building. It is notable that the solving time for both models is fairly comparable, as is shown in Figure 9. This emphasizes the impact of strong mathematical interconnectedness, as in the self-sufficient building model, on the overall solving time, whereas large-scale but easy-to-decompose energy system models such as the dispatch model are relatively quickly solvable compared to their overall size, but consume a large amount of memory.

In addition, a comparison between the aggregation using *typical days* and *typical time steps* reveals that the memory consumption for a certain number of total time steps is, in the case of the dispatch model, almost identical for both the *typical time steps* and *typical days*. For instance, the boxplots for 40 *typical days* and 960 *typical time steps* in the lower part of Figure 10 are almost identical. This highlights that the assumption that the complexity of the model scales with the total number of time steps is valid if no storage components are taken into account. In the case of the self-sufficient building, however, an overhead due to more complex coupling storage conditions can be observed in the case of *typical time steps* compared to the application of *typical days*. Combined with the high deviations induced by an aggregation to *typical time steps* and a long solving time, this approach appears to be highly inconvenient for models with storage components. Moreover, the storage conditions for coupling aggregated *typical time steps* lead to a constant offset to the complexity of an energy system model.

In summary, the findings imply that the overall memory consumption of an energy system model during the solving process is composed of a summand that scales with the number of considered time steps, and a summand that depends on the model's intrinsic number of coupling constraints, e.g., for modeling states of charge of storage components over the entire time horizon, as well as an additional offset introduced by other processes.

# 5. Summary

The results outlined in the preceding section reveal that the clustering of time series to *typical time steps* consistently leads to better clustering indicators in *a priori* input data analyses of the clustered time series compared to clustered *typical days*. However, the analysis of the results reveal that *typical time steps* are not necessarily the superior aggregation method. Specifically, we answer the initially posed research questions as follows:

***Research Question 1:*** *Can an optimal choice of time-series aggregation technique be made a priori based on its capacity to represent the original time series?*

We demonstrate that this is not generally possible: although *typical time steps* outperformed *typical days* in the *a priori* input data analysis for both investigated models, this did not lead to more accurate results for the self-sufficient building model. However, smaller clustering indicators for a predefined model and period length can be used to indicate smaller deviations in the optimization results.

***Research Question 2:*** Are aggregation techniques predetermined for certain types of energy system optimization models and underlying research questions?

With respect to the scope of our analysis, we can state that temporal aggregation techniques based on *typical time steps* outperform *typical days* on the temporally-decoupled dispatch model, whereas the application of *typical days* is the pareto-optimal method for the self-sufficient building model, which takes a variety of different storage systems into account. This is even more remarkable as both *typical time steps and typical days* are modeled in such a way that storage capacities can be considered [81].

These findings are summarized in Figure 11, i.e., the choice of whether to use typical time steps or typical periods depends on the role of storage components or other temporally-coupled constraints and clustering indicators should only be used for comparing different clustering configurations if the model and period lengths are fixed.

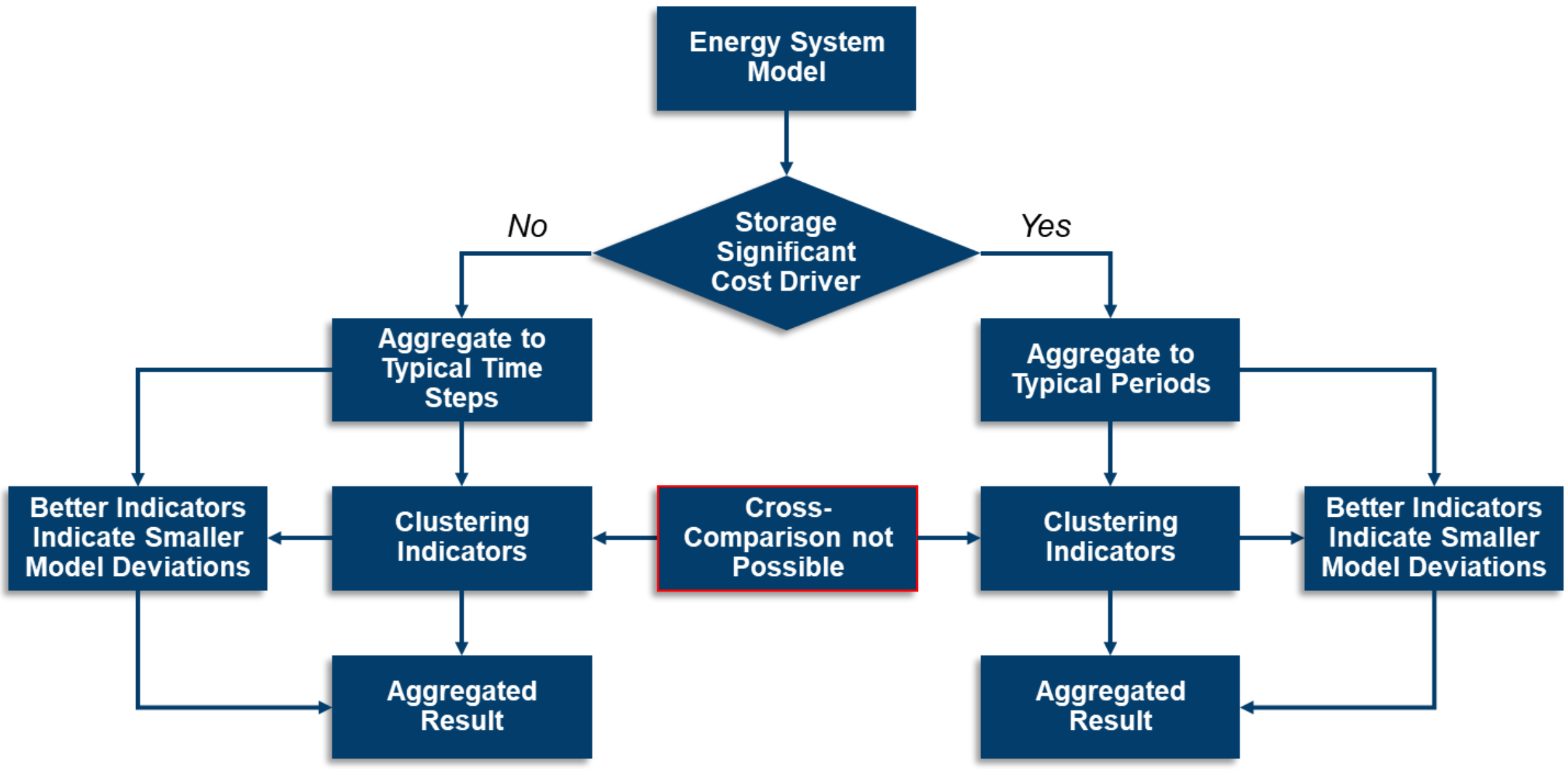

***Figure 11. Summary of the findings.***

## 6. Conclusions

In this study, we analyzed the impact of different temporal aggregation techniques on different energy system models with real-world applications. More specifically, we aggregated the time series to *typical time steps* and, alternatively, to *typical days*. We then fed the reduced time series into a cogeneration system model for a self-sufficient building and a multi-regional electricity dispatch model, and tested the effects on accuracy and complexity.

The findings of this work place common assumptions regarding temporal aggregation for energy system models into a new context. First, the autocorrelation of time series, e.g., a daily pattern of solar and electricity profiles, is not of major importance for the question to what length time periods should be aggregated for an energy system model. The more important indicator is the existence of storage components and the duration of their cycles. Second, it was clearly shown that only a loose relationship exists between clustering indicators and the quality of the aggregated energy system's optimal solution.

Hence, we conclude that knowledge of the mathematical structure of the energy system optimization is necessary to make optimal choices with regard to the temporal aggregation method. In contrast to existing streamlining efforts in the literature that strive to define optimal temporal aggregation methods based on certain features of the time series, our results indicate that the structure of the optimization problem must be considered in order to make optimal choices.

Future research could therefore strive to define the optimal time series optimization methods based on features of the input data time series and the model's specific mathematical features.

## Acknowledgement

This work was supported by the Federal Ministry for Economic Affairs and Energy of Germany in the project METIS (project number 03ET4064A).